\documentclass[12pt,twoside]{article} 
\usepackage{amsfonts}
\usepackage{graphicx}
\usepackage[small,sc]{caption2} 

\newtheorem{theorem}{Theorem}[section]
\newtheorem{lemma}[theorem]{Lemma}
\newtheorem{proposition}[theorem]{Proposition}

\newtheorem{definition}[theorem]{Definition}
\newtheorem{conjecture}[theorem]{Conjecture}
\newtheorem{rmrk}[theorem]{Remark}

\newcommand{\R} {{\mathbb R}}

\newcommand{\Z} {{\mathbb Z}}
\newcommand{\N} {{\mathbb N}}
\newcommand{\qed} {\hfill {\small Q.E.D.} \par\medskip}
\newcommand{\skippar} {\par\medskip}

\newcommand{\proof} {\noindent \textsc{Proof.} }
\newcommand{\proofof}[1] {\noindent \textsc{Proof of {#1}.} }
\newcommand{\article}[3] {\textsc{{#1}}, {\itshape {#2}}, {{#3}}.}
\newcommand{\book}[3] {\textsc{{#1}}, {\itshape {#2}}, {{#3}}.}
\newcommand{\vol} {\textbf}
\newcommand{\eps} {\varepsilon}
\newcommand{\rset}[2] {\left\{ #1 \: \left| \: #2 \right. \! \right\} }
\newcommand{\lset}[2] {\left\{ \left. \! #1 \: \right| \: #2 \right\} }

\renewcommand{\iff} {if and only if\ }

\newcommand{\fn} {function}
\newcommand{\bi} {billiard}
\newcommand{\me} {measure}
\newcommand{\tr} {trajector}
\newcommand{\erg} {ergodic}
\newcommand{\sy} {system}
\newcommand{\hyp} {hyperbolic}
\newcommand{\sca} {scatterer}
\renewcommand{\o} {orbit}

\newcommand{\si} {\mathcal{S}}	
\newcommand{\Sc} {\mathcal{O}}	
\renewcommand{\a} {\alpha}	
\renewcommand{\b} {\beta}	
\newcommand{\is} {\mathcal{I}}	
\newcommand{\ps} {\mathcal{M}}	
\newcommand{\nps} {\mathcal{N}}	
\newcommand{\x} {x}		
\newcommand{\y} {y}		
\newcommand{\ts} {\mathrm{T}}	
\newcommand{\ph} {\varphi}	
\newcommand{\wu} {W^u}		
\newcommand{\ws} {W^s}		
\newcommand{\wsu} {W^{s(u)}}	
\newcommand{\g} {\gamma}	
\newcommand{\clg} {\mathcal{L}}	
\newcommand{\rlg} {\mathcal{R}}	
\renewcommand{\l} {\ell}	
\newcommand{\bfe} {\mathbf{E}}
\newcommand{\bfn} {\mathbf{N}}
\newcommand{\bfw} {\mathbf{W}}
\newcommand{\bfs} {\mathbf{S}}
\newcommand{\bfr} {\mathbf{R}}
\newcommand{\bff} {\mathbf{F}}
\newcommand{\bfl} {\mathbf{L}}
\newcommand{\bfb} {\mathbf{B}}

\newcommand{\sect}[1] {\section{{#1}} \setcounter{equation}{0}}

\newcommand{\fig}[3] {
\medskip\smallskip
\begin{figure}[ht]
	\centering
	\includegraphics[width=#2]{#1.eps}
	\begin{minipage}[t]{0.80\linewidth} 
		\caption{#3}
		\protect\label{#1}
	\end{minipage}
\end{figure}
\medskip
}

\newenvironment{remark}
{\begin{rmrk} \em}
{\end{rmrk}}

\begin{document}

\title{\textbf{Typicality of recurrence for Lorentz gases}}

\author{\textsc{Marco Lenci}
\thanks{
Department of Mathematical Sciences,
Stevens Institute of Technology, 
Hoboken, NJ 07030, U.S.A. \ 
E-mail: \texttt{mlenci@math.stevens.edu} 
} }

\date{October 2004}

\maketitle

\begin{abstract}
	It is a safe conjecture that most (not necessarily periodic)
	two-dimensional Lorentz gases with finite horizon are
	recurrent. Here we formalize this conjecture by means of a
	stochastic ensemble of Lorentz gases, in which i.i.d.\ random
	scatterers are placed in each cell of a co-compact lattice in
	the plane.

	We prove that the typical Lorentz gas, in the sense of Baire,
	is recurrent, and give results in the direction of showing
	that recurrence is an almost sure property (including a
	zero-one law that holds in every dimension).  A few toy models
	illustrate the extent of these results.

\end{abstract}

\sect{Introduction}
\label{sec-intro}

A Lorentz gas (LG) is the \bi\ \sy\ in the complement of a union of
disjoint bounded, regular, convex sets of the plane. Namely, a
dimensionless particle moves with constant unit velocity until it hits
one of the sets (henceforth `\sca s'), at which point it
undergoes an instantaneous Fesnel reflection, i.e., the angle of
reflection equals the angle of incidence.

This dynamical \sy\ generalizes on the one hand the so-called
`Sinai \bi', in which the particle is confined to a bounded
domain, and on the other hand the \emph{periodic} Lorentz gas, in
which the \sca\ configuration is invariant for the action of a
co-compact lattice in $\R^2$.

The most important feature that comes with a LG being an extended \sy\
is that its physically relevant \me, the Liouville \me, is infinite.
In this paper we are interested in the most fundamental \erg\ property
of an infinite-\me\ \sy: Poincar\'e recurrence. This property is far
from trivial, in general \cite{a}.  As a matter of fact, it took the
community considerable effort and more than a decade to prove
recurrence for periodic LGs with finite horizon (i.e., such that the
free path between two collisions is bounded from above). This was
achieved by Schmidt \cite{sch} and Conze \cite{co} at the end of the
1990's. Recurrence has earned further importance lately, as it was
proved \cite{l} that it is a sufficient condition for \erg ity. (See
Section \ref{sec-defs} for the definition of \erg ity in infinite \me\
and for the general geometric assumptions that are needed for
this and all the forthcoming results.) It is well known that, for
\emph{dispersing} \bi s such as ours, \erg ity implies much stronger
chaotic properties---in our case, for suitable finite-\me\ Poincar\'e
maps. (Good surveys of old and recent results in this celebrated field
are found in \cite{ks}, \cite{sc}, \cite{cy}, \cite{cm}).

\skippar

The motivation behind the present work is the idea that ``most''
finite-horizon Lorentz gases must be recurrent. After all, if the most
orderly \sca\ configurations, the periodic ones, give rise to
diffusive and thus recurrent dynamical \sy s, one imagines that the
same must happen for the ``typical'' configuration. This conviction is
corroborated by the results of \cite{l}. For instance, a LG can be
non-recurrent only if it is totally transient, i.e., almost all \tr
ies escape to infinity. Also, a compactly supported perturbation of a
recurrent LG is recurrent as well. As a matter of fact, no example has
been constructed yet of a transient LG.

In order to formalize the intuition above into a precise conjecture, we
present a very natural space of LGs, an \emph{ensemble}, in the sense
that it comes endowed with a probability \me. Given a lattice with
compact fundamental domain, we partition the plane into copies of this
domain (henceforth `cells'). In each cell we place a random
configuration of \sca s so that the configurations in two distinct
cells are independent and identically distributed. We call this
ensemble $\clg$.

The conjecture then reads: almost all gases in $\clg$ are
recurrent. We are not able to prove the conjecture as of now, but can
give a list of results that will hopefully put it within closer reach
(and anyway make it all the more credible). For instance, the set
$\rlg$ of recurrent LGs is topologically typical in $\clg$, provided
one metrizes the latter in a reasonable way (cf.\ Section
\ref{sec-model}). Moreover, $\rlg$ has either full or zero \me\
(Section \ref{sec-meas}). This last result is particularly valuable as
it generalizes to all dimensions.

We also construct a finite-\me\ dynamical \sy\ whose dynamics
comprises that of all \o s in all configurations of $\clg$. The almost
sure recurrence in $\clg$ is equivalent to the \emph{cocycle}
recurrence of a certain \fn\ over this \sy, and we give a sufficient
condition for that (Section \ref{sec-finite}). Verifying the condition
on our model, however, seems rather complicated, so we study how this
dynamical \sy\ behaves for a few simple models (Section
\ref{sec-toy}).

From a technical viewpoint, the paper builds on the results of
\cite{l}. The reader, however, need not know the details of the
proofs, but just the statements, which are given for convenience in
Section \ref{sec-defs}.

\bigskip

\noindent
\textbf{Acknowledgments.} I wish to thank Lai-Sang Young, Fran\c cois
Ledrappier and Charles Newman for very useful discussions. This work
was partially supported by NSF Grant DMS-0405439. Previous travel
funding from GNFM (Italy) is also acknowledged.

\sect{Definitions and preliminary results}
\label{sec-defs}

Let $\{ \Sc_\a \}_{\a\in\is}$ be a family of pairwise disjoint, open,
bounded, convex subsets of $\R^2$, with $C^3$ boundary; $\is$, the
index set, is assumed countable.

With the term `Lorentz gas' (LG) we will indicate both the family $\{
\Sc_\a \}$ (also called `\sca\ configuration' or simply
`configuration') and the \bi\ \sy\ in $\R^2 \setminus
\bigcup_{\a\in\is} \Sc_\a$. The following definitions, assumptions and
basic facts regarding the \bi\ dynamics are standard---and, at any
rate, described in larger detail in \cite{l}---therefore we will lay
them out rather concisely.

To each \sca\ one associates the cylinder $\ps_\a := S^1_{L_\a} \times
[0,\pi]$, where $S^1_{L_\a}$ is the circle of circumference $L_\a$,
the latter being the length of $\partial \Sc_\a$. A pair $(r,\ph) \in
\ps_\a$ represents the element $(q,v)$ of the unit tangent bundle of
$\R^2 \setminus \bigcup_\a \Sc_\a$ thus determined: $q$ is the point
of $\partial \Sc_\a$ parametrized by the arc-length coordinate $r$ (an
origin $r=0$ is fixed once and for all on every $\Sc_\a$, and $r$
increases when moving counterclockwise along $\partial \Sc_\a$); $v
\in \ts_q \R^2$ is the unit vector based in $q$ that forms a
counterclockwise angle $\ph$ with the tangent line to $\Sc_\a$ at $q$,
and points outwardly w.r.t.\ $\Sc_\a$ (see Fig.~\ref{ftyp1}). In the
rest of the paper will also denote pairs $(r,\ph)$ by $x$.

\fig{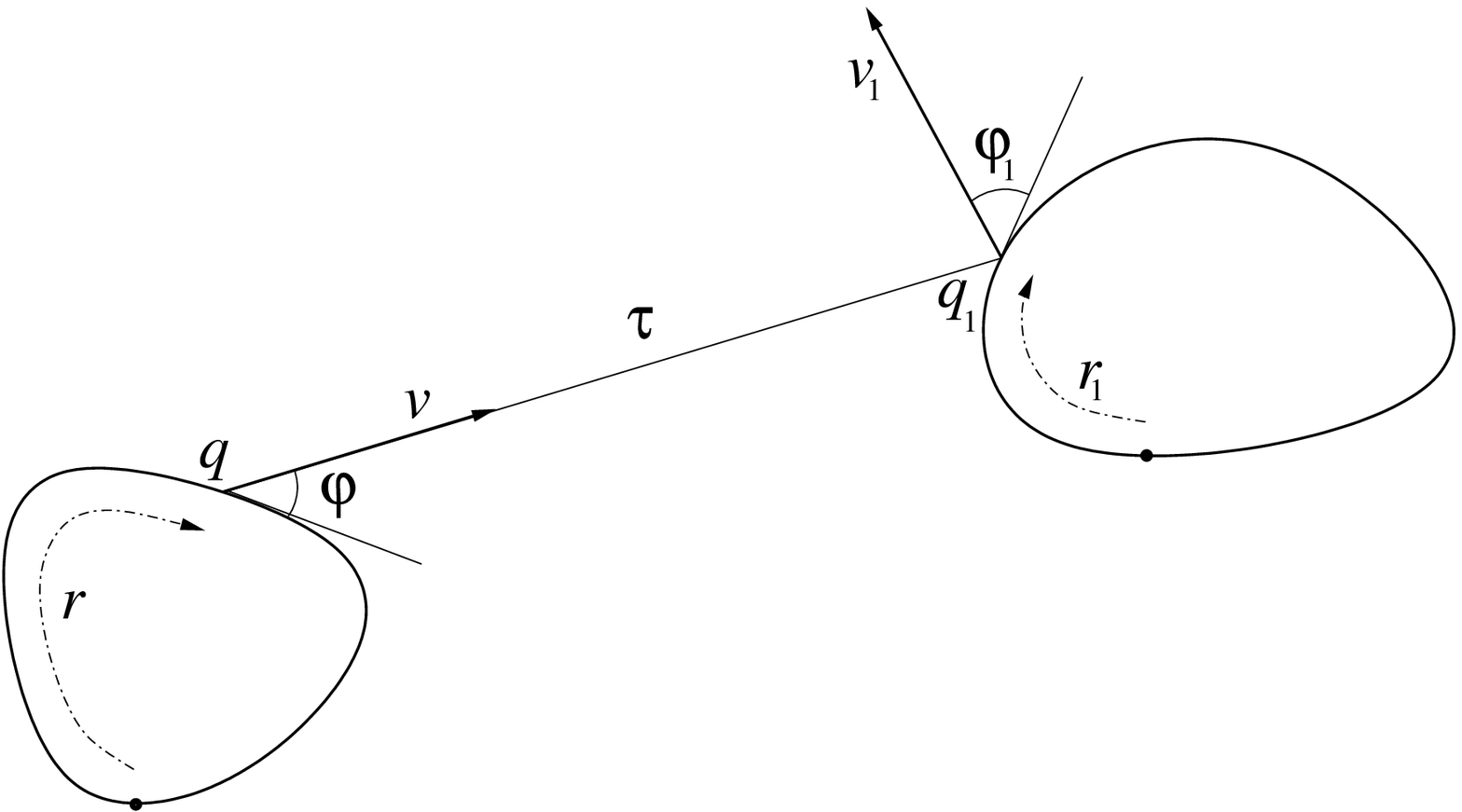} {4.3in} {Basic definitions for the \bi\ map.}

The phase space $\ps := \bigsqcup_\a \ps$ is the disjoint union of all
the $\ps_\a$ (the disjoint union is a needed formality because points
belonging to different $\ps_\a$ may be denoted by the same pair
$(r,\ph)$). We introduce a map $T$ whose action is illustrated in
Fig.~\ref{ftyp1}: $T \x = \x_1$ if $\x$ represents $(q,v)$, $\x_1$
represents $(q_1,v_1)$, and a material point in $q$ traveling with
velocity $v$ has its first collision at $q_1$ with postcollisional
velocity $v_1$. $T$, which is called the \bi\ map, preserves a \me\
$\mu$ on $\ps$, defined by the density $d\mu(r,\ph) := \sin\ph \, dr
d\ph$. Clearly $\mu(\ps) = 2 \sum_\a L_\a = \infty$, save for some
pathological situations when the size of $\Sc_\a$ accumulates at
zero. As a matter of fact, such situations will be explicitly excluded
by the following assumptions, that we maintain throughout the paper.

If $\x$ represents $(q,v)$ with $q \in \ps_\a$, let $k(\x)$ denote the
curvature of $\partial \Sc_\a$ in $q$, and $\tau(\x)$ the \emph{free
path} of $\x$, i.e., the distance between $q$ and $q_1$, the next
collision point (see again Fig.~\ref{ftyp1}). By hypothesis, there
exist $k_m, k_M, \tau_m, \tau_M >0$, such that, $\forall \x \in \ps$,
\begin{eqnarray}
	&& k_m \le k(\x) \le k_M;		\label{cond-k} \\ 
	&& \tau_m \le \tau(\x) \le \tau_M.	\label{cond-tau}
\end{eqnarray}
The second inequality in (\ref{cond-tau}) is the celebrated
\emph{finite horizon} condition. It is clear that (\ref{cond-k})
implies that the size of any $\Sc_\a$ is bounded above and below.

The following definitions may not be obvious for dynamical \sy s of
infinite \me:

\begin{definition}
	The \me-preserving dynamical \sy\ $(\ps, T, \mu)$ is called
	\linebreak[4]
	\textbf{(Poincar\'e) recurrent} if, for every measurable $A
	\subseteq \ps$, the \o\ of $\mu$-almost every $\x \in A$
	returns to $A$ at least once (and thus infinitely many times,
	due to the invariance of $\mu$).  
	\label{def-rec}
\end{definition}

\begin{definition}
	The \me-preserving dynamical \sy\ $(\ps, T, \mu)$ is called
	\textbf{\erg} if every $A \subseteq \ps$ measurable and
	invariant mod $\mu$ (i.e., $\mu( T^{-1} A \,\Delta\, A) = 0$),
	has either zero \me\ or full \me\ (i.e., $\mu (\ps \setminus
	A) = 0$).  
	\label{def-erg}
\end{definition}

Dispersing \bi s like the \sy\ at hand are prototypical examples of
\hyp\ \sy s with singularities. The presence of the singularities
represents a conspicuous hurdle in proving the \hyp\ and \erg\
properties. This is even more so in infinite \me\ and, as is the case
here, when the singularities themselves have an infinite extension (in
the sense of their length as smooth curves in $\ps$).

The following three results are the technical backbone of \cite{l}:

\begin{theorem}
	The Lorentz gas introduced above has a \textbf{\hyp\
	structure}, meaning that for $\mu$-a.e.\ $\x \in \ps$ there
	are local stable and unstable manifolds (LSUMs) at $\x$,
	denoted $\wsu(\x)$. These two measurable foliations (when
	endowed with a Lebesgue-equivalent transversal \me) are
	absolutely continuous w.r.t.~$\mu$.  
	\label{thm-hyp}
\end{theorem}

For a precise definition of LSUM in this context, see \cite{l}. Here
we are primarily interested in their core property: if $\y \in
\wsu(\x)$ then $d_\ps(T^n \x, T^n \y) \to 0$ as $n \to
+\infty(-\infty)$; $d_\ps$ is the Riemannian distance in $\ps$
(by definition $d_\ps (\x,\y) = \infty$ if $\x$ and $\y$ belong to
different cylinders $\ps_\a$, $\ps_\b$). For our \sy s one can see
that the rate of vanishing is exponential.

\begin{theorem}
	Given $\a \in \is$, almost every two points $\x,\y \in \ps_a$
	are connected by a polyline of alternating LSUMs, in the sense
	that there is a finite collection of LSUMs $\ws(\x_1)$,
	$\wu(\x_2), \ws(\x_3), \,\cdots , \wu(\x_m)$, with $\x_1 :=
	\x$ and $\x_m := \y$, such that each LSUM intersects the next
	transversally.  
	\label{thm-conn}
\end{theorem}

\begin{theorem}
	$(\ps, T, \mu)$ is \erg\ \iff it is recurrent.
	\label{thm-rec-erg}
\end{theorem}

For this last theorem, only the sufficient condition was given in
\cite{l}, but the necessary condition is obvious, anyway: if there is
a positive-\me\ wandering set, one can split it in two non-trivial
parts, which must necessarily belong to two different \erg\
components.  

\skippar

The following statement was used (and justified) in \cite{l}, but
never explicitly emphasized.

\begin{proposition}
	A LG as introduced above is either recurrent, i.e., totally
	conservative, or totally dissipative.  
	\label{prop-rec-diss}
\end{proposition}

\proof The dissipative part $D$ of $(\ps, T, \mu)$ is defined as the
maximal countable union of wandering sets of $\ps$, modulo $\mu$
\cite{a}. If $\mu(D)>0$, we claim that $D$ contains whole LSUMs; that
is, $\ws(D) := \rset{y \in \ws(\x)} {\x \in D}$ and the analog
$\wu(D)$ are equal to $D$, modulo $\mu$. We prove that first
statement, the second being obviously equivalent.

Take a positive-\me\ wandering set $A$. Without loss of generality, $A
\subseteq \ps_\a$ for some $\a \in \is$. Apart from a null-\me\ set,
$A$ is the disjoint union of
\begin{equation}
	A_n := \lset{\x \in A} {n = \max \lset{k \ge 0} {T^k \x \in
	\ps_\a}},
\end{equation}
with $n \ge 0$ (it is easy to see that almost no points of $A$ can
return to $\ps_\a$ infinitely many times). Pick $n$ for which
$\mu(A_n)>0$. $\ws(T^n A_n)$ is a wandering set, because points in the
same LSM have the same forward itinerary w.r.t.\ the partition $\{
\ps_\b \}$ (that is, they hit the same \sca s in the future); in
particular, if $\y \in \ws(\x)$ with $\x \in T^n A_n$, then $T^k \y
\not\in \ps_\a$, for all $k>0$.

As is known, the local stable foliation can be chosen invariant, that
is, $T \ws(\x) \subseteq \ws(T \x)$---this is in fact a standard
assumption.  Therefore $T^n \ws(A_n) \subseteq \ws(T^n A_n)$. Together
with the above conclusions, this implies that $T^n \ws(A_n)$, and thus
$\ws(A_n)$, is wandering. Repeating the argument for all $n$ such that
$\mu(A_n)>0$ proves that $\ws(A)$ is wandering, yielding our initial
claim.

By Theorem \ref{thm-conn}, then, any $\ps_\a$ is either wholly
contained in $D$ or in its complement. If $D \ne \ps$, there must be
two nearest neighbors $\Sc_\a$ and $\Sc_\b$ such that $\ps_\a
\subseteq D$ and $\ps_\b \cap D = \emptyset$. But this is absurd as
$D$ is $T$-invariant and there exists $B \subset \ps_\a$, with
$\mu(B)>0$, such that $T B \subset \ps_\b$.  
\qed

In this paper we are interested in recurrent LGs, so let us start to
give examples thereof. Recall the definition of periodic LG from the
Introduction.

\begin{theorem}
	\emph{\cite{sch, co}} A periodic LG with finite horizon and
	strictly convex \sca s is recurrent.
	\label{thm-per}
\end{theorem}

\begin{definition}
	The LG $\{ \Sc_\a \}_{\a \in \is}$ is called a \textbf{finite
	modification} of $\{ \Sc_\a \}_{\a \in \is_0}$ if $\is =
	(\is_0 \setminus \is_1 ) \cup \is_2$, where:
	\begin{itemize} 
		\item[(a)] $\is_1$ is a finite subset of $\is_0$.

		\item[(b)] $\is_2$ is the index set of a finite LG
		such that $d_{\R^2} (\Sc_\a, \Sc_\b) > 0$ for any $\a
		\in \is_2$, $\b \in \is_0 \setminus \is_1$ ($d_{\R^2}$ 
		is the distance in the plane).
	\end{itemize}
	\label{def-fin-mod}
\end{definition}

\begin{proposition}
	A LG is recurrent \iff any of its finite modifications are
	recurrent. 
	\label{prop-fin-mod}
\end{proposition}

\proof The proof of the necessary condition is identical to that of
Proposition 5.3 of \cite{l}.  The sufficient condition follows
automatically since a LG is a finite modification of any of its finite
modifications.
\qed

\sect{Model and topological typicality}
\label{sec-model}

As explained in the Introduction, it is reasonable to conjecture that
most LGs are recurrent. We need a satisfyingly general class of gases
for which `most' can be properly defined. Our choice is explained
hereafter.

Consider a co-compact lattice $\Gamma \subset \R^2$, with $\{ C_\g
\}_{\g \in \Gamma}$ its corresponding partition of the plane, that is,
$\R^2 = \bigcup_{\g \in \Gamma} C_\g$, with $C_\g = C_0 + \g$, and
$C_\g \cap C_\eta = \emptyset$ for $\g \ne \eta$. In each \emph{cell}
$C_\g$ we put a random configuration of \sca s parametrized by
$\l_\g$, where $\{ \l_\g \}$ are independent identically distributed
random variables from the same probability space $(\Omega, \pi)$. (In
the remainder, a generic element of $\Omega$ with be denoted by
$\omega$.)  We assume that (\ref{cond-k})-(\ref{cond-tau}) are
satisfied for every realization of this random field.

Examples are illustrated by Fig.~\ref{ftyp2} and its caption.

\fig{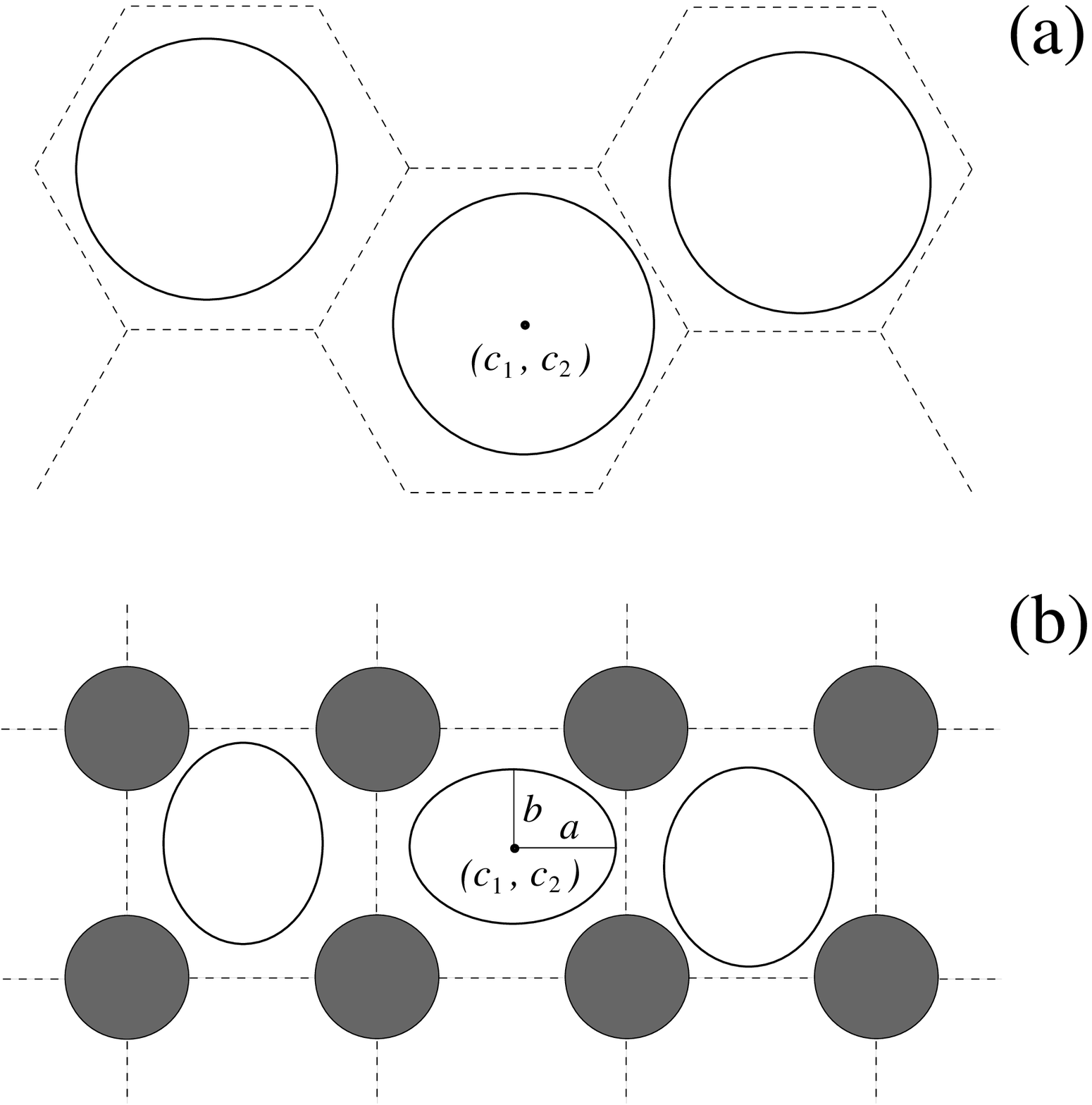} {4.1in} {Two examples of random LGs.  In (a), $\Gamma =
\mathrm{Hex}$; in each cell the \sca\ is a disc of radius $R$ and
random center $(c_1, c_2) =: \omega \in \Omega := B(0, r)$, with $r$
sufficiently small; $\pi$ is the normalized Lebesgue \me\ on $\Omega$.
In (b), $\Gamma = \Z^2$; the \sca\ is an ellipse of random center
$(c_1, c_2)$ and random semiaxes $a,b$, with $\omega := (c_1, c_2, a,
b) \in B(0, r) \times I_1 \times I_2 =: \Omega$ ($I_1, I_2$ are
intervals); $\pi$ is the normalized Lebesgue \me\ on $\Omega$; the
gray, non-random \sca s are needed to comply with the finite-horizon
condition.}

The class of LGs we will concern ourselves with is $(\clg, \Pi) :=
(\Omega, \pi)^\Gamma$, where the superscript denotes the product of
$\Gamma$ copies of $(\Omega, \pi)$. From now on a Lorentz gas will be
an element $\l = \{ \l_\g \} \in \clg$.

\skippar

In many cases, just as in the examples of Fig.~\ref{ftyp2}, $\Omega$
is also a metric space. We ask that the metric verifies the following
natural property.

\begin{definition}
	A distance \fn\ $d_\Omega$ on pairs of $\Omega$ is called
	\textbf{compatible with the dynamics} if:
	\begin{itemize} 
		\item[(a)] $(\Omega,d_\Omega)$ is a compact metric 
		space.

		\item[(b)] Every \sca\ $\Sc^{(i)} (\omega)$ $(i = 1,
		\ldots , N)$ represented by $\omega \in \Omega$
		depends in a $C^3$ fashion on $\omega$. In other
		words, if $C \subset \R^2$ is a cell and
		$\xi_\omega^{(i)}: S^1 \longrightarrow C$ is the
		arc-length parametrization of $\partial \Sc^{(i)}
		(\omega)$, renormalized to 1, then $\|
		\xi_\omega^{(i)} - \xi_{\omega'}^{(i)} \|_{C^3(S^1)}
		\to 0$, when $d_\Omega( \omega, \omega') \to 0$.
	\end{itemize} 
	\label{def-comp}
\end{definition}

Obviously, if $\Omega$ is finite, then $d_\Omega$ is compatible with
the dynamics.

\skippar

The above definition induces a distance on $\clg$ which makes it a
complete metric space. Namely, for $\l = \{ \l_\g \}$ and $\l' = \{
\l'_\g \}$,
\begin{equation}
	d_\clg (\l, \l') = \sum_{\g \in \Gamma} 2^{-|\g|} \, d_\Omega
	(\l_\g, \l'_\g) 
	\label{d-clg}
\end{equation}

In this setup, recurrence is a typical property, in the sense of
Baire:

\begin{theorem}
	If $d_\Omega$ is compatible with the dynamics, then
	\begin{displaymath} 
		\rlg := \lset{\l \in \clg} {\l \mbox{ \rm is a
		recurrent LG}} 
	\end{displaymath} 
	contains a $G_\delta$-set.
	\label{thm-top-typ}
\end{theorem}

\proof We simplify the proof somewhat if we consider the cylinder
\begin{equation}
	\clg_0 := \rset{\l = \{ \l_\g \} } {\l_0 = \omega_0}, 
\end{equation}
where $\omega_0$ is a fixed element of $\Omega$. In view of
Proposition \ref{prop-fin-mod}, Theorem \ref{thm-top-typ} is
equivalent to showing that $\rlg_0 := \rlg \cap \clg_0$ is a residual
set of $\clg_0$ in the appropriate topology.

First of all, let us construct a countable subset of $\rlg_0$ that is
dense in $\clg_0$. Take a dense sequence $\{ \omega_j \}_{j \in \N}$
in $\Omega$ (if $\Omega$ is finite, one can use $\Omega$ instead of
$\{ \omega_j \}$). Let $\Lambda$ denote a finite subset of $\Gamma
\setminus \{ 0 \}$, and $\{ j_\eta \}_{\eta \in \Lambda} \in
\N^\Lambda$ an $m$-tuple of natural numbers indexed by the elements of
$\Lambda$ (here $m = \# \Lambda$). To each pair $(\Lambda, \{ j_\eta
\}) =: n$ is associated the configuration $\l^{(n)}$, defined by
\begin{equation}
	\l_\g^{(n)} = \left\{
	\begin{array}{rcl}
		\omega_{j_\g}, && \mbox{if } \g \in \Lambda; \\
		\omega_0, && \mbox{if } \g \not\in \Lambda.
	\end{array}
	\right.
\end{equation}
Since the set of such pairs $n$ is countable, let us pretend that $n
\in \N$. By looking at definition (\ref{d-clg}), it is rather clear
that $\{ \l^{(n)} \}_{n \in \N}$ is dense in $\clg$. Furthermore, each
$\l_n$ is a finite modification of a periodic LG; hence $\l^{(n)} \in
\rlg_0$.

Now, for any $\l$, let us consider a specific \sca\ $\Sc_0$ in the
cell $C_0$; for instance, $\Sc_0 = \Sc^{(1)} (\omega_0)$ (which, in
the notation of Definition \ref{def-comp}\emph{(b)}, means the
``first'' \sca\ of $C_0$). The crucial point is that $\Sc_0$ is
exactly the same for every $\l \in \clg_0$, because the configuration
in $C_0$ is fixed. We naturally call $\ps_0$ the cylinder in phase
space corresponding to $\Sc_0$. Also, set $\mu_0 ( \,\cdot\, ) := \mu
( \,\cdot\, ) / \mu(\ps_0)$.

If $\mathcal{B}(\l, \rho) \subset \clg_0$ denotes the ball of center
$\l$ and radius $\rho>0$, w.r.t.\ $d_\clg$, we contend that for
every $n,m \in \N$ there exists $\eps_m^{(n)}>0$ such that, for
all $\l \in \mathcal{B}(\l^{(n)}, \eps_m^{(n)})$ the set
\begin{equation}
	A(\l) := \rset{\x \in \ps_0} {\exists k>0 \mbox{ such that }
	T_\l^k x \in \ps_0} 
	\label{def-a-l}
\end{equation}
has \me\ 
\begin{equation}
	\mu_0 \left( A(\l) \right) \ge \left( 1 - \frac1m \right). 
	\label{me-a-l}
\end{equation}
(In (\ref{def-a-l}), $T_\l$ represents the \bi\ map for the LG $\l$.)

It is not too hard to verify this claim, once we have unraveled its
rather intricated formulation. In fact, $\l^{(n)}$ is recurrent and
thus $A(\l^{(n)})$ has full \me\ in $\ps_0$. Take then $\x \in
A(\l^{(n)})$, with $k>0$ its first return time to $\ps_0$. The \tr y
of $x$ up to $T_{\l^{(n)}}^k \x$ is non-singular in the sense that its
polyline representation on $\R^2 \setminus \bigcup_\a \Sc_\a$ is
tangent to no $\Sc_\a$ (by convention, singular \tr ies, a null-\me\
set, are ignored, at least after they hit the tangency). Therefore,
if one slightly modifies the shape and location of the scatterers of
$\l^{(n)}$ (thus turning it into some $\l$ with $d_\clg(\l, \l^{(n)})
< \eps$), then the sequence of \sca s hit by $\{ T_\l^j \x \}_{j=0}^k$
is the same as for $\{ T_{\l^{(n)}}^j \x \}_{j=0}^k$. In particular
$T_\l^k \x \in \ps_0$. (For the \emph{cogniscenti}: $\x \in \ps_0$
will have the same forward itinerary up to time $k$, w.r.t.\ the
partition $\{ \ps_\a \}$, if the perturbation of the LG, which induces
a perturbation on the singularity set $\si$, leaves $\x$ in the same
connected component of $\ps_0 \setminus \si_k$, where $\si_k := \si
\cap T^{-1} \si \cap \,\cdots\, \cap T^{-k+1}$.  Sufficiently small
perturbations of $\l^{(n)}$ will obviously do this.)
 
The above reasoning shows that, for all $\x \in A(\l^{(n)})$, there
exists $\eps = \eps(\x)>0$ such that, $\forall \l \in \mathcal{B}
(\l^{(n)}, \eps)$, $\x \in A(\l)$, too. Whence the claim.

Finally, the set
\begin{equation}
	\mathcal{G} := \bigcap_{m\in\N} \, \bigcup_{n\in\N} \,
	\mathcal{B}(\l^{(n)}, \eps_m^{(n)})
\end{equation}
is $G_\delta$ by construction. From (\ref{me-a-l}), $\mu_0 (A(\l)) =
1$, for all $\l \in \mathcal{G}$. This proves Theorem
\ref{thm-top-typ} because, for such $\l$, it follows that $\ps_0$
belongs in the conservative part of $(\ps_\l, T_\l, \mu)$ (with the
obvious meaning for $\ps_\l$). Therefore, by Proposition
\ref{prop-rec-diss}, $\l \in \rlg_0$.  
\qed

\begin{remark}
	It is evident that the lattice structure of $\clg$ played
	essentially no role in the proof of Theorem
	\ref{thm-top-typ}. Using the same method, one can prove
	the same result for any complete metric space $\mathcal{X}$ of
	LGs such that:
	\begin{itemize}
		\item $\mathcal{X}$ has a dense set of recurrent gases
		(e.g., finite modifications of periodic LGs).

		\item The distance is compatible with the dynamics,
		that is, two configurations $\l, \l' \in \mathcal{X}$
		are close \iff there is a bijective correspondence
		between their \sca s such that corresponding \sca s
		have $C^3$-uniformly close boundaries.
	\end{itemize}
\end{remark}

\sect{Measure-theoretic typicality}
\label{sec-meas}

We are aiming for a stronger notion of typicality for the recurrence
property, as the space $\clg$ was constructed with a built-in
probability \me\ $\Pi$.

\begin{conjecture}
	$\Pi(\rlg) = 1$.
	\label{conj}
\end{conjecture}

This seems very credible, especially in light of Theorem
\ref{thm-per}: if a periodic configuration produces a recurrent
dynamics, then a typical random configuration will randomize the
motion of the particle even more, making it possibly even more similar
to a random walk.

Unfortunately, Conjecture \ref{conj} will remain such throughout the
paper. The following result, however, seems to indicate that we are on
the right track.

\begin{theorem}
	If $\mathcal{A}$ is the $\sigma$-algebra induced on $\clg$ by
	its construction (i.e., $\mathcal{A} = \mathcal{C}^{\otimes
	\Gamma}$, where $\mathcal{C}$ is the $\sigma$-algebra defined
	on $\Omega$), then $\rlg \in \overline{\mathcal{A}}^\Pi$ and
	$\Pi(\rlg) \in \{0,1\}$.
	\label{thm-0-1}
\end{theorem}

\proof The second assertion is rather trivial once we establish the
first. In fact, consider this natural action of $\Gamma$ on $\clg$:
for $\eta \in \Gamma$,
\begin{equation}
	\sigma_\eta (\l) =: \l' = \{ \l'_\g \}_{\g \in \Gamma}, \quad
	\mbox{with} \quad \l'_\g := \l_{\g+\eta}.
	\label{def-sigma}
\end{equation}
Obviously, $\sigma$ preserves the \me\ $\Pi$. Furthermore,
$\sigma_\eta(\rlg) = \rlg$ for all $\eta \in \Gamma$, since
recurrence is a translation invariant property. On the other hand,
$(\clg, \{ \sigma_\eta \}_{\eta \in \Gamma}, \Pi)$ is \erg\
(it is by definition a generalized Bernoulli shift in two
dimensions). These two facts imply that $\Pi(\rlg) \in \{0,1\}$.

\skippar

For the first statement we use the same trickery as in the proof of
Theorem \ref{thm-top-typ}. From Proposition \ref{prop-fin-mod} we know
that $\rlg$ is invariant w.r.t.\ changes in the $0^\mathrm{th}$
component (i.e., $\l \in \rlg \, \Longleftrightarrow \l' \in \rlg$,
for all $\l'$ such that $\l'_\g = \l_\g$, whenever $\g \ne
0$). More in detail, it is a ``cylinder'' whose sections are the
$\rlg_0$ introduced in the proof of Theorem \ref{thm-top-typ} (one for
each $\omega_0$). If we call $\mathcal{A}_0$ and $\Pi_0$,
respectively, the factor $\sigma$-algebra and the factor \me\ induced
by $\mathcal{A}$ and $\Pi$ on the cylinder $\clg_0$ (notice that
$(\clg_0, \mathcal{A}_0, \Pi_0) \simeq (\Omega, \mathcal{C},
\pi)^{\Gamma \setminus \{ 0 \}}$), then
\begin{equation}
	\rlg_0 \in \overline{\mathcal{A}_0}^{\Pi_0} \quad
	\Longrightarrow \quad \rlg \in \overline{\mathcal{A}}^\Pi.
\end{equation}
As for proving the above l.h.s., we recall the definition of $\ps_0$
from the proof of Theorem \ref{thm-top-typ}, and set
\begin{equation}
	A := \rset{(\x,\l) \in \ps_0 \times \clg_0} {\limsup_{k \to
	+\infty} \: (\chi_{\ps_0} \circ T_\l^k) (\x) =1 }.
\end{equation}
$A$ is measurable because $T_\l\, \x$ is clearly a measurable \fn\ of
$(\x, \l)$ (indeed, due to the finite-horizon condition, it does not
depend on the \sca s of $\l$ that are at a certain distance from
$\Sc_0$; so it is even measurable w.r.t.\ a certain subalgebra of
sets depending only on a finite number of lattice sites).

Now, Proposition \ref{prop-rec-diss} implies that, for any given $\l$,
either a full-\me\ or a zero-\me\ set of points in $\ps_0$ come back
to $\ps_0$ infinitely many times, depending on $\l$ being recurrent or
not. This amounts to saying that, almost surely, $A$ contains whole
``horizontal'' fibers of $\ps_0 \times \clg_0$, that is, $A = \ps_0
\times \rlg_0$ mod $\mu \times \Pi$. By Lemma \ref{lemma-prod} of
the Appendix, $\rlg_0 \in \overline{\mathcal{A}_0}^{\Pi_0}$.  
\qed

\begin{remark}
	Theorem \ref{thm-0-1} is much more general than was presented
	here, and applies easily to the $d$-dimensional case. In fact,
	the only non-trivial ingredient in the proof is Proposition
	\ref{prop-rec-diss}, which is in turn a consequence of Theorem
	\ref{thm-conn} (that is just a weak formulation of the local
	ergodicity theorem). Therefore, if (\ref{cond-k}) is
	substituted by
	\begin{equation}
		k_m \le \mathbf{k}(q) \le k_M,
		\label{ddim-cond-k}
	\end{equation}
	where $\mathbf{k}(q)$ is the second fundamental form of
	$\partial \Sc_\a$ at $q$ (the inequalities here are meant in
	the sense of the quadratic forms), then Theorem \ref{thm-0-1}
	holds for the class $\clg = \clg(d, \Gamma, \Omega, \pi)$ of
	$d$-dimensional LGs with i.i.d.\ random \sca s in every cell
	of $\Gamma$, selected from the probability space $(\Omega,
	\pi)$, whenever the geometry of the \sca s makes the local
	ergodicity theorem hold. This includes at least all
	semi-dispersing \sca s given by algebraic equations
	\cite{bcst}.  Moreover, one can apply this zero-one law to
	many situations in which the dimension of $\Gamma$ is strictly
	less than the dimension of the Euclidean space (e.g., a 3D
	billiard in an infinite parallelepiped acted upon by $\Z$, and
	so on...).  
	\label{rk-0-1}
\end{remark}

\sect{A finite-\me\ dynamical \sy}
\label{sec-finite}

From a technical point of view, the difficulties associated with our
\sy\ arise by and large from the fact that the given invariant \me\
has infinite mass. But the lattice structure of $\clg$ suggests the
construction of a \emph{finite-\me} dynamical \sy\ that embodies
\emph{all} LGs in $\clg$.

Consider the cell $C_0$ associated to the origin of $\Gamma$: we think
of it as our fundamental domain. Call $\partial_{*} C_0$ the part of
$\partial C_0$ that does not intersect any non-random \sca\ ($\partial
C_0$ can never intersect a random \sca, anyway, lest (\ref{cond-tau})
be violated; as a matter of fact, the random \sca s must keep at least
$\tau_m/2$ units away from $\partial C_0$). In example (a) of
Fig.~\ref{ftyp2}, $\partial_{*} C_0 = \partial C_0$, whereas in
example (b), $\partial_{*} C_0$ is the union of four disjoint segments
of equal length. Define
\begin{equation}
	\nps := \rset{(q,v) \in \ts \R^2} {q \in \partial_{*}
	C_0,\ |v|=1, \mbox{ and } v \mbox{ points inwardly w.r.t.} \:
	C_0 }.  
	\label{def-nps}
\end{equation}
To maintain consistency with the notation of Section \ref{sec-defs},
we identify $\partial_{*} C_0$ with a subset $J$ of $\R$, in which an
arc-length coordinate $r$ uniquely determines a point $q \in
\partial_{*} C_0$. Then, if $\ph$ parametrizes the direction of $v$ in
the usual way (like in Fig.~\ref{ftyp1}), then $\nps$ can be
identified with $J \times [0, \pi]$.

\begin{remark}
	This identification is always flawed at a finite number of
	points in $\partial_{*} C_0$. For instance, in
	Fig.~\ref{ftyp2}(a), at the six vertices of $\partial C_0$; in
	Fig.~\ref{ftyp2}(b), at the eight boundary points of
	$\partial_{*} C_0$. There are two ways to do away with this
	problem. The first way is tantamount to ignoring it: one can
	exclude these points from $\partial_{*} C_0$ (in which case,
	$\partial_{*} C_0$ will always be a disjoint union of open
	intervals). This exclusion is acceptable since it affects only
	a null-\me\ subset of $\partial_{*} C_0$, w.r.t.\ the relevant
	\me\ that we introduce below. The second way consists in
	identifying, on a case-by-case basis, different pairs $(r,
	\ph)$ and $(r_1, \ph_1)$, corresponding either to the same
	line element, or to the pre- and post-collisional line
	elements for the same collision. For instance, in example (b),
	if $r_0$ is the left endpoint of an interval of $\partial_{*}
	C_0$, $(r_0, \ph) \simeq (r_0, \pi-\ph)$.
\end{remark}

Let us call $\mu_1$ the standard \bi-invariant \me\ for the
cross-section $\nps$, normalized to 1 (in $(r, \ph)$ coordinates, $d
\mu_1 (r, \ph) = [2 \, \mathrm{length} (\partial_{*} C_0)]^{-1} \sin
\ph \, dr d\ph$). If $\omega \in \Omega$ determines the configuration
of \sca s in $C_0$, we can define a map $R_\omega : \nps
\longrightarrow \nps$ as follows.  Trace the (forward) \tr y of $\x :=
(q,v) \in \nps$ until it crosses $\partial C_0$ for the first
time---see Fig.~\ref{ftyp3}. This occurs at the point $q_1$ and with
velocity $v_1$. Say that $C_\g$ is the cell that the particle enters
upon leaving $C_0$. Define then
\begin{eqnarray}
	R_\omega \, \x = R_\omega (q,v) &:=& (q_1-\g, v_1) \in \nps; 
	\label{def-tomega} \\
	e(\x, \omega) &:=& \g \in G.	\label{def-e}
\end{eqnarray}
Here $G \subset \Gamma$ is the set of \emph{primitive directions} of
$\Gamma$, each corresponding to a neighboring cell of $C_0$. We name
$e$ the `exit \fn'. Finally, $R_\omega$ preserves $\mu_1$. (To
give but a brief explanation, $\partial_{*} C_0$ is a
\emph{transparent wall} for the \bi\ flow.  Poincar\'e maps for
transparent walls are virtually the same as those for reflecting
walls---they are actually a commonly used trick in \bi\ dynamics, cf.\
\cite{l0}.)

\fig{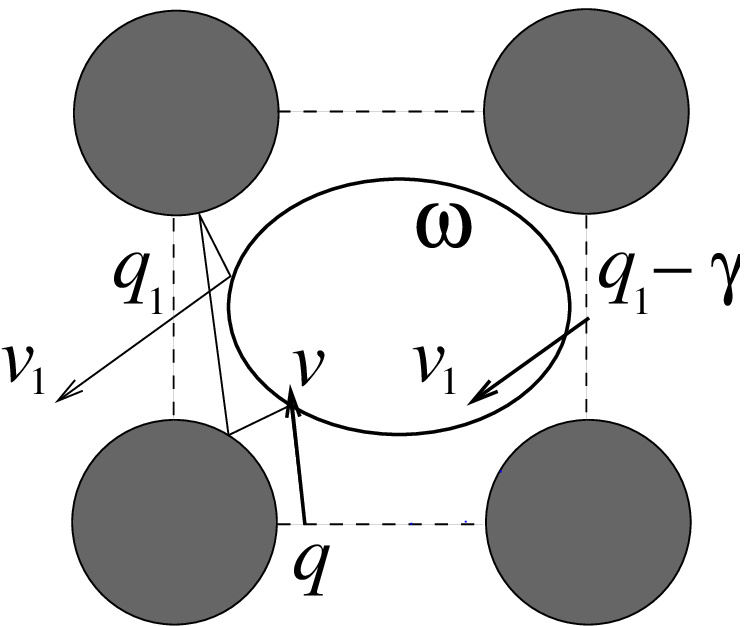} {1.7in} {The definition of $R_\omega$. In this case
$R_\omega (q,v) = (q_1-\g, v_1)$, with $\g = (-1,0) \in \Z^2$.}

The dynamical \sy\ that we want to introduce in this section is the
triple $(\Sigma, F, \nu)$, where:
\begin{itemize}
	\item $\Sigma := \nps \times \clg$.

	\item $F(\x, \l) := \left( R_{\l_0}\, \x, \sigma_{e(\x, \l_0)}
	(\l) \right)$, defining a map $\Sigma \longrightarrow
	\Sigma$. Here $\sigma$ is the $\Gamma$-action on $\clg$
	defined by (\ref{def-sigma}) and $\l_0$ is, as usual, the
	$0^\mathrm{th}$ component of $\l$.

	\item $\nu := \mu_1 \times \Pi$. Since $\mu_1$ is
	$R_\omega$-invariant for every $\omega \in \Omega$, and $\Pi$
	is $\sigma$-invariant, then $\nu$ is $F$-invariant.
\end{itemize}

The idea behind this definition is that, instead of following a given
\o\ form a cell to another, every time we shift the LG in the
direction opposite to the \o\ displacement, so that the point always
lands in $C_0$.  Clearly, $F: \Sigma \longrightarrow \Sigma$
encompasses the dynamics of all points on all LGs of $\clg$.  It is
equally as clear that we are in the case in which a.e.\ $\l \in \clg$
is recurrent \iff the \fn\ $e$ verifies the following:

\begin{definition}
	Let $(\Sigma, F, \nu)$ be a mea\-sure-preserving dynamical
	\sy\ with $\nu(\Sigma) = 1$. If $e: \Sigma \longrightarrow
	\Gamma \subseteq \R^d$, define the \textbf{cocycle}
	\begin{displaymath}
		S_n(z) := \sum_{k=0}^{n-1} (e \circ F^k) (z).
	\end{displaymath}
	The \fn\ $e$ (or the cocycle $S_n$) is called
	\textbf{recurrent} if, for $\nu$-almost all $z \in \Sigma$,
	\begin{displaymath}
		\liminf_{n \to +\infty} \left| S_n(z) \right| = 0.
	\end{displaymath}
	\label{def-f-rec}
\end{definition}

When $\Gamma$ is discrete, which is our case, the above is equivalent
to saying that $S_n(z)=0$ infinitely often in $n$.

A notable sufficient condition for cocycle recurrence was given by
Schmidt:

\begin{theorem}
	\emph{\cite{sch}} Assume that $(\Sigma, F, \nu)$ is \erg, and
	denote by $p_n$ the distribution of $S_n/n^{1/d}$, i.e., for a
	Borel set $A$ of $\R^d$,
	\begin{displaymath} 
		p_n(A) := \nu\left( \rset{z\in \Sigma}{\frac{S_n(z)} 
		{n^{1/d}} \in A } \right).
	\end{displaymath}
	If there exists a positive-density sequence $\{ n_k
	\}_{k\in\N}$ and a constant $c>0$ such that
	\begin{displaymath}
		p_{n_k}(B(0,\rho)) \ge c {\rho^d} 
	\end{displaymath}
	for all sufficiently small balls $B(0,\rho)$ of center 0 and
	radius $\rho$ in $\R^d$, then the cocycle $\{ S_n \}$
	(equivalently, the \fn\ $e$) is recurrent.  
	\label{thm-f-rec}
\end{theorem}

\begin{remark}
	In the case of interest to this paper, that is $d=2$,
	estabilishing the Central Limit Theorem for the family of
	variables $\{ e \circ F^k\}$ (even with a degenerate limit) is
	clearly enough to apply Theorem \ref{thm-f-rec}. This is in
	fact how Schmidt proves Theorem \ref{thm-per} via \cite{bs}.
\end{remark}

Coming back to the actual \sy\ at hand, this is what we know:

\begin{proposition}
	If $(\Sigma, F, \nu)$ is the dynamical \sy\ introduced above
	then
	\begin{itemize}
		\item[(a)] Every measurable invariant set of $\Sigma$
		is of the form $\nps \times B$ mod $\nu$, where $B$ is
		a measurable set of $\clg$. Furthermore, either $B$ or
		$\clg \setminus B$ has empty interior.

		\item[(b)] The \sy\ is topologically transitive.

		\item[(c)] In the case of almost sure recurrence (that
		is, when $\Pi(\rlg) = 1$ or, which is the same, when
		$\{ S_n \}$ is a recurrent cocycle), $(\Sigma, F,
		\nu)$ is \erg.
	\end{itemize}
	\label{prop-f-erg}
\end{proposition}

\proof For a given $\l \in \clg$, consider the dynamical \sy\
$(\ps_\l, T_\l, \mu)$, corresponding to the LG $\l$. In view of
Theorem \ref{thm-hyp}, we construct ``local stable and unstable
manifolds'' for $F$ at a.e.\ point of $\nps \times \{\l\}$. (More
precisely, the ``LUMs'' are constructed as push-forwards of the LSMs
of $T_\l$ onto the cross-section $\nps$; analogously, the ``LUMs'' are
pull-backwards of the LUMs of $T_\l$.) These curves are contained in
$\nps \times \{\l\}$. The quotation marks are in order here as they
are not \emph{bona fide} LSUMs for $(\Sigma, F, \nu)$, which is not a
hyperbolic \sy\ in any reasonable sense.

We now exploit Theorem \ref{thm-conn} to conclude that, in each
connected component of $\nps \times \{\l\}$, a.e.\ pair of points
(w.r.t.\ $\mu_1$) are joined through a polyline of ``LSUMs'',
therefore, via the usual Hopf argument, they lie in the same \erg\
component of $(\Sigma, F, \nu)$. On the other hand, it is easy to
verify that no two connected components of $\nps \times \{\l\}$ (each
corresponding to a different segment of $\partial_{*} C_0$) can have
separate dynamics, i.e., belong to distinct \erg\ components. In
conclusion, at least for a.a.\ $\l \in \clg$, $\nps \times \{\l\}$ is
cointained in the same \erg\ component of $(\Sigma, F, \nu)$.  Which
is to say, the only $F$-invariant sets of $\Sigma$, modulo $\nu$, are
of the form $\nps \times B$. That $B$ is measurable mod $\Pi$ in
$\clg$ is a consequence of Lemma \ref{lemma-prod} of the
Appendix. This proves the first part of statement \emph{(a)}.

\skippar

Next, consider two open sets $U_1, U_2 \in \Sigma$. For the purpose of
proving topological transitivity one can always pass to subsets, so
assume that, for $i\in \{1,2\}$, $U_i = V_i \times B_i$, where $V_i$
is an open set of $\nps$ and $B_i$ is a cylinder of $\clg$. This means
that, for each $i$, there exists a finite subset of $\Gamma$, called
$\Lambda_i$, and a family of open sets of $\Omega$, $\{ A_i^\eta
\}_{\eta \in \Lambda_i}$, such that $B_i = \rset{\l\in\clg} {\forall
\eta \in \Lambda_i, \l_\eta \in A_i^\eta}$.

Take a sufficiently large $\g_0 \in \Gamma$ so that $(\Lambda_2 +
\g_0)$ does not intersect $\Lambda_1$. It is clearly possible to find
a periodic $\bar{\l} = \{ \bar{\l}_\g \}$ such that, for all $\eta \in
\Lambda_1$, $\bar{\l}_\eta \in A_1^\eta$ and, for all $\eta \in
\Lambda_2$, $\bar{\l}_{\eta+\g_0} \in A_2^\eta$. By construction,
$\bar{\l} \in B_1$ and $\sigma_{\g_0} (\bar{\l}) \in B_2$. By Theorem
\ref{thm-per}, $\bar{\l}$ is recurrent so $(\ps_{\bar{\l}},
T_{\bar{\l}}, \mu)$ is \erg. This implies that almost every \bi\ \tr
y, in the three-dimensional phase space of $\l$, intersects the
cross-section defined by $\rset{(q+\gamma_0,v)} {(q,v) \in V_2}$. In
other words, since $\mu_1(V_1)>0$, there exist a non-singular $\bar{x}
\in V_1$ and an integer $n$ such that $F^n(\bar{x}, \bar{\l}) \in V_2
\times \{ \sigma_{\g_0} (\bar{\l}) \} \subset V_2 \times B_2$. But
since $\bar{x}$ is non-singular and the metric on $\Omega$ is
compatible with the dynamics (Definition \ref{def-comp}), we can
perturb $\bar{x}$ and $\bar{\l}$ a little bit and still end up in $V_2
\times B_2$. This means that there exists an open neighborhood
$\mathcal{U}$ of $(\bar{x}, \bar{\l})$ such that $F^n(\mathcal{U})
\subset V_2 \times B_2$. This fact proves \emph{(b)} and the second
part of \emph{(a)}.

\skippar

For \emph{(c)} we use the following lemma, whose proof will be given
below.

\begin{lemma}
	For $\l \in \clg$ and $\g \in \Gamma$, set
	\begin{displaymath}
		D_\l^\g := \rset{\x \in \nps} {S_n(x,\l) = \gamma,\:
		\mbox{\rm for some} \ n \in \N}
	\end{displaymath}
	and
	\begin{displaymath}
		E := \rset{\l \in \clg} {\forall \g\in\Gamma, \: 
		\mu_1(D_\l^\g) > 0}.
	\end{displaymath}
	If $\Pi(E)>0$ then $(\Sigma, F, \nu)$ is \erg.
	\label{lemma-tech}
\end{lemma}

In the hypothesis of \emph{(c)}, a.e.\ $\l$ is an \erg\ Lorentz gas so,
by the argument used earlier, $\mu_1(D_\l^\g)=1$ for all $\g$. Thus
$\Pi(E)=1$.  
\qed

\proofof{Lemma \ref{lemma-tech}} Suppose the \sy\ is not \erg.  By
Proposition \ref{prop-f-erg}\emph{(a)}, we have an invariant set $\nps
\times B$ (mod $\nu$), with $B$ a Borel set of $\clg$ and $\Pi(B) \in
(0,1)$. Set $B^c := \clg \setminus B$.  Either $B$ or $B^c$ (or both)
must intersect $E$ in a positive-\me\ subset. Say that this happens
for $B$. Since $(\clg, \{ \sigma_\eta \}, \Pi)$ is \erg, one can find
$O \subseteq B \cap E$ and $\g \in \Gamma$ such that $\Pi(O)>0$ and
\begin{equation}
	\sigma_\g (O) \subseteq B^c.
	\label{tech1}
\end{equation}
Fix $\l \in O$. The hypotheses of the lemma imply that there is a
positive integer $n$ and a set $D_\l^{\g,n} \subset D_\l^\g$, with
$\mu_1 (D_\l^{\g,n}) >0$, such that $S_n(x,\l) = \g$ for all $x\in
D_\l^{\g,n}$. That is to say,
\begin{equation}
	F^n ( D_\l^{\g,n} \times \{\l\} ) \subseteq \nps \times
	\sigma_\g (\l) \subseteq \nps \times B,
\end{equation}
the last inclusion holding at least for a.e.\ $\l \in O$, due the
$F$-invariance of $\nps \times B$ mod $\nu$ (notice that we have
implicitly used Fubini's Theorem and Proposition
\ref{prop-f-erg}\emph{(a)}). This gives that $\sigma_\g (O) \subseteq
B$ mod $\Pi$, in contradiction with (\ref{tech1}).  
\qed

\sect{Toy models}
\label{sec-toy}

As we have seen, exploring the statistical properties of $(\Sigma, F,
\nu)$ is not exactly a trivial task. In this section we consider much
simplified versions of that dynamical \sy\ that nonetheless have the
same lattice structure. One can call this structure `deterministic
dynamics in a random enviroment'. The intent is to get an idea of
those properties of the \sy\ that depend more on the random
environment than on the details of the dynamics. (A less easy model
will be treated in \cite{ln}.)

In these examples $\Gamma$ will always be $\Z^2$. Let us denote
\begin{equation}
	\bfe = g_1 = (1,0), \quad \bfn = g_2 = (0,1), \quad \bfw = g_3
	= (-1,0), \quad \bfs = g_4 = (0,-1),
	\label{enws}
\end{equation}
the symbols standing for East, North, West and South. These are the
primitive directions of $\Z^2$ and together they form the set $G$. To
each of these directions is associated a copy of the unit square
$[0,1]^2$. These four copies are named $\nps_\bfe = \nps_1$,
$\nps_\bfn = \nps_2$, $\nps_\bfw = \nps_3$, and $\nps_\bfn = \nps_4$;
also $\nps := \bigsqcup_{i=1}^4 \nps_i$. A point $x \in \nps_\bfe$
corresponds to the particle entering $C_0$ from the western side (its
incoming direction being $\bfe$), and so on analogously. We endow
$\nps$ with $\mu_1$, the Lebesgue \me\ divided by 4, which is the
right normalization factor here.

To complete the definition of our toy version of $(\Sigma, F, \nu)$,
following the paradigm of Section \ref{sec-finite}, we need to
introduce the probability space $(\Omega, \pi)$ that governs the
randomness of each cell; the map $R_\omega: \nps \longrightarrow \nps$
that gives the dynamics in a cell in the state $\omega \in \Omega$;
and the exit \fn\ $e: \nps \times \Omega \longrightarrow G$.  All
these objects will vary from example to example.

Before discussing the models one by one, let us introduce the
$m$-baker's map $K_m: [0,1]^2 \longrightarrow [0,1]^2$.  For $m$ a
positive integer and $y := (y_1,y_2) \in [0,1]^2$,
\begin{equation}
	K_m(y_1,y_2) = \left( \{ m y_1 \} , \frac{y_2 + [m y_1]} m
	\right),
\end{equation}
where $[\rho]$ and $\{\rho\}$ are the integer and fractional part,
respectively, of $\rho \ge 0$. Of course, $K_2$ is the standard
baker's map.

\subsection{Example 1: Mimicking the standard random walk}
\label{subs-ex1}

Here $R_\omega$ acts on each $\nps_i$ as a $4$-baker's map. Precisely,
if $x := (y,i) \in \nps$ (with $y \in [0,1]^2$ and $i \in \{1, \ldots
, 4\}$), then
\begin{equation}
	R_\omega(y,i) := \left( K_4(y), e((y,i),\omega) \right).
\end{equation}
Thus the dynamics does not really depend on the random state of the
cell: $R$ depends on $\omega$ only through $e$, that is, only insofar
as $R_\omega(x)$ must necessarily land on $\nps_{e(x,\omega)}$.  

In this example, $\Omega := \{ 1,2,3 \}$ and $\pi$ is any probability
\me\ there (there is nothing special about the number 3, and that is
exactly the point in choosing it).  The exit \fn\ $e$ is given in
terms of level sets by Fig.~\ref{ftyp4}.  Observe that $e$ depends
non-trivially on $\omega$, otherwise the \sy\ would be of a much
simpler nature and recurrence would just amount to the \fn\ recurrence
of $e$ relative to $(\nps,R,\mu_1)$.

\fig{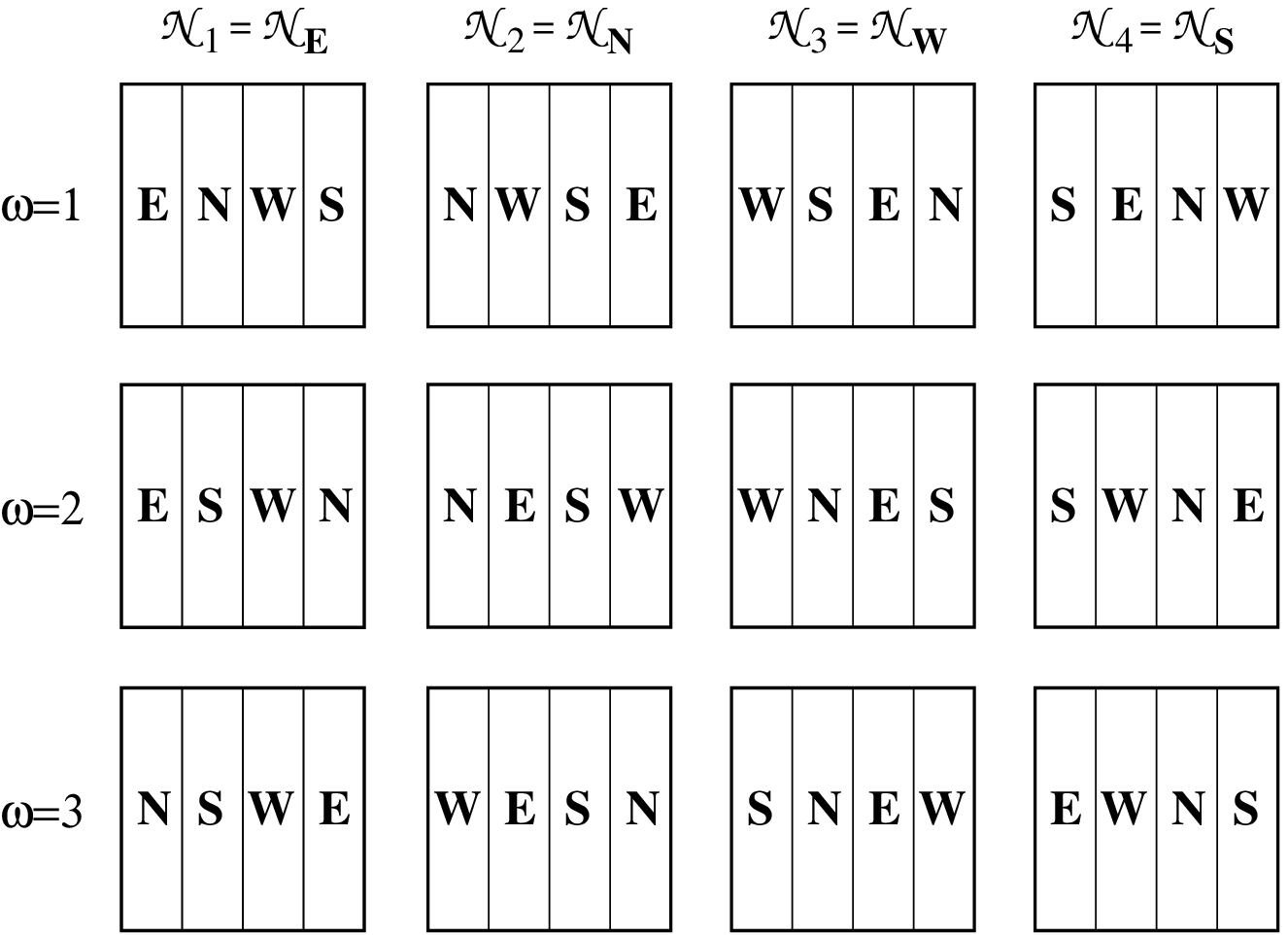} {4.6in} {The definition of $e$ for Example 1. Each row
displays a copy of $\nps := \bigsqcup_{i=1}^4 \nps_i$, corresponding
to different values of $\omega \in \{ 1,2,3 \}$. The level sets are
all rectangles of base $1/4$ and height 1. The label in each level set
is the common image of the set via $e$, according to the
notation (\ref{enws}).}

\begin{remark}
	Notice that, for all $\omega \in \Omega$ and $i = 1, \ldots,
	4$, $R_\omega$ is invertible and 
	\begin{equation}
		\mu_1( \rset{x\in\nps} {e(x,\omega) = g_i} ) =
		\frac14. 
	\end{equation}
	These facts are essential if we want to think of our dynamical
	\sy\ as generated by a ``physical'' (read: conservative and
	invertible) \sy\ $(\ps, T, \mu)$, as described in Section
	\ref{sec-finite}. The same will occur also for Examples 2 and
	3.
\end{remark}

We claim that in this setup the particle moves as in a standard random
walk for every realization of the random enviroment. In more exact
terms, for every $\l \in \{ 1,2,3 \}^{\Z^2}$, the stochastic process
\begin{equation}
	(\nps, \mu_1) \ni x \longmapsto \left\{ S_n(x,\l)
	\right\}_{n\ge 0} = \left\{ \sum_{k=0}^{n-1} (e \circ F^k)
	(x,\l) \right\}_{n\ge 0}
\end{equation}
is a standard random walk. In fact, if $\{ \g_k \}_{k=0}^n$ is a path
in $\Z^2$ (i.e., $\g_0=0$ and $|\g_{k+1} - \g_k| = 1 \ \forall k$)
and $i \in \{ 1, \ldots, 4 \}$, one realizes that the conditional
probability
\begin{equation}
	\mu_1 \left( \left. (e \circ F^n) (\,\cdot\, , \l) = g_i \,
	\right| \, S_k(\,\cdot\, , \l) = \g_k, \ \forall k=1, \ldots,
	n \right) = \frac14.
\end{equation}
(The set in which we condition is a rectangle of base $4^{-n+1}$ and
height 1. The preimages of $e \circ F^n$ will subdivide it into 4
rectangles of base $4^{-n}$ and height 1, one for each $g_i$.)

Recurrence is thus guaranteed in \emph{every} fiber $\nps \times \{ \l
\}$, which is an even stronger statement than we sought.

\subsection{Example 2: Left-Right random walk}
\label{sub-ex2}

The previous example was indeed much too easy, and we did not utilize
at all the considerations of Section \ref{sec-finite}. Example 2 is
going to be a tad more involved.  Here $\Omega = \{ 1,2 \}$ and, once
again, the choice of $\pi$ is irrelevant; $e$ is given by
Fig.~\ref{ftyp5}.  $R_\omega$ acts as the standard baker's map, and
its dependence on $\omega$ is as trivial as in the previous model:
\begin{equation}
	R_\omega(y,i) := \left( K_2(y), e((y,i),\omega) \right).
\end{equation}

\fig{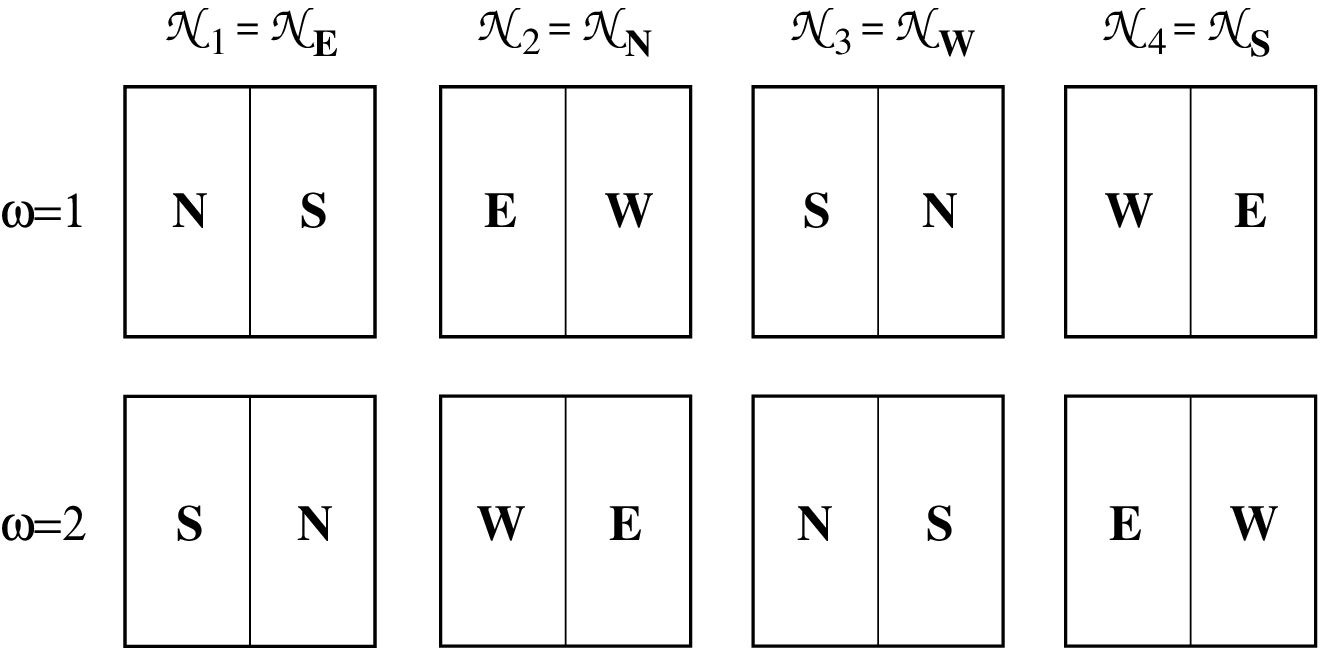} {4.6in} {The definition of $e$ for Example 2. See caption
of Fig.~\ref{ftyp4} for explanations.}

Reasoning along the same lines as in Section \ref{subs-ex1}, we see
that, fixing $\l \in \Omega^{\Z^2}$ and letting $x$ range randomly in
$\nps$ according to $\mu_1$, we obtain the so-called Left-Right random
walk. This is the stochastic process in a which a particle moves in
$\Z^2$ turning its direction by 90 degrees after every step, with a
fifty-fifty chance of turning left or right (left and right being
relative to the direction of the motion; the absolute directions in
$\Z^2$ are $\bfw, \bfe$ and so on).

What makes this model more complicated than Example 1 is that it is
not a Markov chain (at least not in its simplest formulation, that is,
as a random process whose $n^\mathrm{th}$ component is the position of
the particle at time $n$).  Not that probabilists have a hard time
proving the recurrence of this model, but here we will do so by
applying Theorem \ref{thm-f-rec}.

First of all, we claim that $(\Sigma, F, \nu)$ is \erg. This might not
be so apparent, as the dynamics has an obvious symmetry. In fact,
denote $\nps^{(1)} := \nps_\bfe \sqcup \nps_\bfw$ and $\nps^{(2)} :=
\nps_\bfn \sqcup \nps_\bfs$. For a typical $\l \in \clg$, the
evolution of $\nps^{(1)} \times \{\l\}$ is always separated from
$\nps^{(2)} \times \{\l\}$ because, if $x \in \nps^{(1)}$, it is easy
to check from Fig.~\ref{ftyp5} that $F^n (x, \l)$ belongs to
$\nps^{(1)} \times \clg$ or $\nps^{(2)} \times \clg$ depending on $n$
begin even or odd, respectively; but in order for the point to come
back to the same cell (i.e., $F^n (x, \l) \in \nps \times \{\l\}$) $n$
must be even, and therefore $F^n (x, \l) \in \nps^{(1)} \times
\{\l\}$.

This notwithstanding, that $(\Sigma, F, \nu)$ is \erg\ can be seen as
follows: For $j = 1,2$, define $\Sigma^{(i)} := \nps^{(j)} \times
\clg$ and consider the dynamical \sy\ $(\Sigma^{(j)}, F^2, \nu/2)$. It
is easy to ascertain that the ``horizontal'' fibers $\nps^{(j)} \times
\{\l\}$ are wholly contained in the \erg\ components of the \sy. Then
we can use Lemma \ref{lemma-tech} with
\begin{equation}
	\Gamma = \Z_{\mathrm{even}}^2 := \rset{\g = (\g_1, \g_2) \in
	\Z^2} {\g_1 + \g_2 \in 2\Z}
\end{equation}
to conclude that $(\Sigma^{(j)}, F^2, \nu/2)$ is \erg. (The hypothesis
of the lemma applies as it is clear that, for all $\l \in \clg$, one
can reach any $\g \in \Z_{\mathrm{even}}^2$ for some---thus
many---initial conditions $x \in \nps ^{(j)}$. The proof works because
$(\clg, \{ \sigma_\eta \}_{\eta \in \Gamma}, \Pi)$ is \erg\ for
$\Gamma = \Z_{\mathrm{even}}^2$ as well.)  Now, if $U$ is an
$F$-invariant subset of $\Sigma$ mod $\nu$, set $U^{(j)} := U \cap
\Sigma^{(j)}$. Clearly $F U^{(1)} = U^{(2)}$ and $F^2 U^{(1)} =
U^{(1)}$ (mod $\nu$). Hence, both $U^{(j)}$ have either \me\ zero or
full \me\ in $\Sigma^{(j)}$, whence the \erg ity of $(\Sigma, F,
\nu)$.

Now, for every $\l$, it is obvious that the projections of $S_n(
\,\cdot\, , \l)$ onto the horizontal and vertical directions of $\Z^2$
are independent one-dimensional random walks (at times $[n/2]$ or
$[n/2]+1$), which verify the 1D Central Limit Theorem for $n \to
+\infty$. Their orthogonal sum must then verify the 2D Central Limit
Theorem, which implies the same for the ``more random'' process $S_n(
\,\cdot\, , \,\cdot\,)$. This, together with the \erg ity of $(\Sigma,
F, \nu)$, gives the recurrence via Theorem \ref{thm-f-rec},

\subsection{Example 3: Deterministic walk in a random environment}

The next and last example shows that recurrence in $\clg$ is not due
solely to the chaotic nature of the dynamics (which tends to produce a
diffusive behavior in \emph{every} LG), but may also be a consequence
of the random environment.  To demonstrate this point, which makes
Conjecture \ref{conj} all the more convincing, we take $R$ to be as
regular as it can be, practically the identity. Let us define
\begin{equation}
	R_\omega(y,i) := \left(y, e(i,\omega) \right);
\end{equation}
that is, $y$ remains constant and plays absolutely no role, not even
on the exit \fn\ $e$, which, for a given state $\omega$ of the cell,
depends only on the incoming direction $i$.  For all practical
purposes, then, each $\nps_i$ can be collapsed to a point. We actually
do so and for the remainder of the section we consider $\Sigma := \{
1,2,3,4 \}$ or, equivalently, $\Sigma := \{ \bfe, \bfn, \bfw, \bfs
\}$. Instead of a `dynamics in a random environment', we have a
`\emph{walk} in a random environment'.

Let $\Omega := \{ 1,2,3,4 \}$ and, for $\omega \in \Omega$,
\begin{equation}
	e(i,\omega) := i + \omega \mbox{ mod } 4,
\end{equation}
where `mod 4' means the congruent integer between 1 and 4. If we
rename the elements of $\Omega$ as $\bfl := 1$, $\bfb := 2$, $\bfr :=
3$ and $\bff := 4$ (the symbols standing for Forward, Backward, Left
and Right), we see that every time the particle reaches a cell in the
state $\omega$, it will take a step in the direction indicated by
$\omega$ (relative to the incoming direction). Let us denote
$\pi_\bfl, \pi_\bfb, \pi_\bfr, \pi_\bff$ the probabilities of each of
the four symbols.

It is clear that a walk is recurrent \iff it is a periodic \o\ of
$(\Sigma, F, \nu)$. It is also clear that Theorem \ref{thm-0-1} does
not hold---or better, its proof does not apply: recurrence here is not
a translation invariant property. The same for Proposition
\ref{prop-f-erg}\emph{(a)} and the analog of Proposition
\ref{prop-rec-diss}: the existence of a closed walk tells us nothing
about the other walks. We collect what we know in a proposition.

\begin{proposition} 
	Using the notation of the previous sections on the \sy\
	$(\Sigma, F, \nu)$ defined above,
	\begin{itemize}
		\item[(a)] $\Pi(\rlg)>0$.

		\item[(b)] There exists a number $p_c \in (1/2, 1)$
		such that, if $\pi_\bfb > p_c$, or $\pi_\bfl > p_c$, 
		or $\pi_\bfr > p_c$, then $\Pi(\rlg)=1$.

		\item[(c)] If $\pi_\bfl + \pi_\bfr =1$ (and thus
		$\pi_\bfb = \pi_\bff = 0$), then $\Pi(\rlg)=1$.
	\end{itemize}
\end{proposition}

\proof Part \emph{(a)} is obvious since, for every initial direction
$i$, one can always fix $\l_\g$ for a finite number of cells near the
origin so as to create a periodic \o. The LGs $\l$ with those
components fixed form a positive-\me\ cylinder.

As for \emph{(b)}, $p_c$ is the critical probability for the site
percolation in $\Z^2$ (believed to be approximately 0.59). If the
probability of a given $\omega$ is bigger than this number, there is
almost surely a closed loop of cells marked $\omega$ that surrounds
the origin \cite{g}. If $\omega \in \bfb, \bfl, \bfr$, this prevents
any path starting at the origin from reaching infinity (cf.\
\cite{bt}).

Statement \emph{(c)} is Theorem 2(ii) of \cite{bt}.
\qed

\appendix

\sect{Appendix: A lemma from measure theory}

\begin{lemma}
	Let $(X, \mathcal{A}, \mu)$ and $(Y, \mathcal{B}, \nu)$ be two
	probability spaces, and let $\mathcal{A} \otimes \mathcal{B}$
	denote the $\sigma$-algebra on $X \times Y$ generated by the
	rectangles (i.e., sets of the type $A \times B$, with $A \in
	\mathcal{A}$, $B \in \mathcal{B}$).

	If $A \in \mathcal{A}$, with $\mu(A)>0$, and $A \times B \in
	\mathcal{A} \otimes \mathcal{B}$, then there exists a $B_0 \in
	\mathcal{B}$ such that $B \,\Delta\, B_0$ is a subset of a
	$\nu$-null-\me\ set. In particular, if $(Y, \mathcal{B}, \nu)$
	is complete, then $B \in \mathcal{B}$.
	\label{lemma-prod}
\end{lemma}

\proofof{Lemma \ref{lemma-prod}} We first notice that, for every $C
\in \mathcal{A} \otimes \mathcal{B}$, there exists a sequence $\{ U_n
\}$, where $U_n$ is a finite union of rectangles, such that $\lim_{n
\to +\infty} (\mu \times \nu) (C \,\Delta\, U_n) =0$.

We prove the above by showing that the class $\mathcal{C}$ of such
sets $C$ that can be approximated by finite unions of rectangles is
a $\sigma$-algebra in $\mathcal{A} \otimes \mathcal{B}$. Since
$\mathcal{C}$ contains the rectangles, it must be the whole
$\mathcal{A} \otimes \mathcal{B}$.

Obviously $\mathcal{C}$ is closed by complementation (an
approximating sequence for $C^c$ being $U_n^c$, also a finite unione
of rectangles) and by countable disjoint union (it suffices to neglect
the sets in the ``tail'' of the countable union).

So, let $U_n$ be an approximating sequence for $A \times B$. For $y
\in Y$, denote by
\begin{equation}
	S_{y,n} := \rset{x \in A} {(x,y) \in U_n}
\end{equation}
the section of $U_n \cap (A \times Y)$ relative to $y$. Set $m_n(y) :=
\mu(S_{y,n})$. This is clearly a measurable \fn\ $Y \longrightarrow
\R_0^{+}$ (it is actually a simple \fn). Since it is also bounded,
then $m_n \in L^1(Y,\nu)$.

Assume for the moment that $x \in A$, $y \in B$. Then $(x,y) \in (A
\times B) \,\Delta\, U_n \Longleftrightarrow (x,y) \not\in U_n
\Longleftrightarrow x \not\in S_{y,n}$. Now assume that $x
\in A$, but $y \not\in B$. By the same token $(x,y) \in (A \times B)
\,\Delta\, U_n \Longleftrightarrow (x,y) \in U_n \Longleftrightarrow x
\in S_{y,n}$. These considerations, combined with Fubini's Theorem,
give
\begin{eqnarray}
	&& (\mu \times \nu) \, [ ((A \times B) \,\Delta\, U_n) \cap (A
	\times Y) ] = \nonumber \\
	&=& \int_B \mu(A \setminus S_{y,n}) \, d\nu(y) + \int_{B^c}
	\mu(S_{y,n}) \, d\nu(y) = \\
	&=& \int_Y |m_n(y) - \mu(A) \, \chi_B(y)| \,
	d\nu(y). \nonumber
\end{eqnarray}
Therefore $m_n \to \mu(A) \chi_B$ in $L^1(Y,\nu)$, as $n \to
+\infty$. For a subsquence, then---let us call it again $\{ m_n
\}$---the convergence occurs almost everywhere. More precisely, there
exists a set $Y_0$, with $\nu(Y_0)=1$, such that, for every $y \in
Y_0$, 
\begin{equation}
	\lim_{n \to +\infty} m_n (y) = \mu(A) \, \chi_B(y).
	\label{lemma-prod-10}
\end{equation}
For a fixed $\rho \in (0, \mu(A))$, we denote by $L_n := \lset{y \in
Y} {m_n(y) \ge \rho}$ a certain filled level set of $m_n$. Then
(\ref{lemma-prod-10}) implies that
\begin{equation}
	Y_0 \cap \bigcup_{m\in\N} \, \bigcap_{n\ge m} L_n = Y_0 \cap B.
\end{equation}
Setting $B_0 := Y_0 \cap B \in \mathcal{B}$ concludes the proof of the
lemma.  
\qed

\end{document}